\documentclass[10pt]{amsart}
\usepackage{latexsym, amsmath,amssymb}

\usepackage{xcolor}

\newcommand{\Rnu}{{\mathcal R}_\nu}

\renewcommand{\epsilon}{\varepsilon}
\newtheorem{theorem}{Theorem}
\newtheorem{lemma}[theorem]{Lemma}
\newcommand{\bcor}{\begin{corr}}

\newcommand{\bdeff}{\begin{deff}}

\newcommand{\bprop}{\begin{proposition}}
\newcommand{\ele}{\end{lemma}}
\newcommand{\ecor}{\end{corr}}
\newcommand{\edeff}{\end{deff}}

\newcommand{\eprop}{\end{proposition}}

\newcommand{\Rn}{{\mathbb R}^n}

\newcommand{\la}{\lambda}

\newcommand{\e}{\varepsilon}

\renewcommand{\Pi}{\varPi}

\renewcommand{\epsilon}{\varepsilon}

\newcommand{\R}{{\mathbb R}}

\begin{document}

\title[]{Approximating pointwise products of\\ Laplacian eigenfunctions}
\thanks{J.L.~is supported in part by the NSF (DMS-1454939);
  C.D.S.~is supported in part by the NSF (NSF Grant DMS-1665373)
  and the Simons Foundation; S.S.~is supported in part by the NSF
  (DMS-1763179) and the Alfred P. Sloan Foundation.}

\author{Jianfeng Lu}
\address[J. L.]{Department of Mathematics, Department of Physics, and Department of Chemistry,
Duke University, Box 90320, Durham NC 27708, USA}
\email{jianfeng@math.duke.edu}

\author{Christopher D. Sogge} \address[C.S.]{Department of
  Mathematics, Johns Hopkins University, Baltimore, MD 21218, USA}
\email{sogge@jhu.edu}

\author{Stefan Steinerberger}
\address[S.S.]{Department of Mathematics, Yale University, New Haven, CT 06511, USA}
\email{stefan.steinerberger@yale.edu}


\begin{abstract} We consider Laplacian eigenfunctions 
  on a $d-$dimensional bounded domain $M$ (or a $d-$dimensional
  compact manifold $M$) with Dirichlet conditions. These
  operators give rise to a sequence of eigenfunctions
  $(e_\ell)_{\ell \in \mathbb{N}}$. We study the subspace of all
  pointwise products
$$ A_n = \mbox{span} \left\{ e_i(x) e_j(x): 1 \leq i,j \leq n\right\} \subseteq L^2(M).$$
Clearly, that vector space has dimension $\mbox{dim}(A_n) = n(n+1)/2$. We prove that products $e_i e_j$ of eigenfunctions are simple in a certain sense: for any $\varepsilon > 0$, there exists a low-dimensional vector space $B_n$ that almost contains all products. More precisely, denoting the orthogonal projection $\Pi_{B_n}:L^2(M) \rightarrow B_n$, we have
$$ \forall~1 \leq i,j \leq n~  \qquad \|e_ie_j - \Pi_{B_n}( e_i e_j) \|_{L^2} \leq \varepsilon$$
and the size of the space $\mbox{dim}(B_n)$ is relatively small: for every $\delta > 0$,
$$ \mbox{dim}(B_n) \lesssim_{M,\delta}  \varepsilon^{-\delta} n^{1+\delta}.$$
We obtain the same sort of bounds for products of arbitrary length, as
well for approximation in $H^{-1}$ norm.  Pointwise products of
eigenfunctions are low-rank. This has implications, among other
things, for the validity of fast algorithms in electronic structure
computations.

\end{abstract}

\maketitle

\section{Introduction}

Let $(M,g)$ be a compact Riemannian manifold of dimension $d\ge 2$.  Let $\{e_j\}$ be an orthonormal basis of eigenfunctions with
frequencies $\la_j$ arranged so that $0=\la_0<\la_1\le \la_2\le \cdots$.  Thus, 
$$-\Delta_g e_j=\la_j^2 e_j, \qquad \langle e_j, \, e_k\rangle= \int_M e_j \, \overline{e_k}\, dV_g=\delta^k_j.$$
Here $\Delta_g$ denotes the Laplace-Beltrami operator and $dV_g$ is the volume element associated with the metric $g$ on $M$.  For simplicity,
we are only considering eigenfunctions of Laplace-Beltrami operators but, using the same proof, all of our results extend to eigenfunctions of second order elliptic operators
which are self-adjoint with respect to a smooth density on $M$.
As we shall see in the final section we can also handle Dirichlet eigenfunctions in domains.\\

We ask a very simple question: what can be said about the function $e_{i}(x) e_{j}(x)$?
Clearly, by $L^2-$orthogonality, the function
$e_{i}(x) e_{j}(x)$ has mean value 0 if $i\neq j$ but what else
can be said about its spectral resolution, for example, the size of
$\left\langle e_i e_j, e_k\right\rangle$?  There are very few
results overall; some results have been obtained in the presence of
additional structure assumptions on $\Omega$ connected to number
theory (see Bernstein \& Reznikoff \cite{bern}, Kr\"otz \& Stanton
\cite{krotz} and Sarnak \cite{Sarnak}). Already the simpler question
of understanding $L^2-$size of the product is highly nontrivial: a
seminal result of Burq-G\'{e}rard-Tzetkov \cite{burq} states
$$ \|e_{\mu}e_{\lambda} \|_{L^2} \lesssim \min( \lambda^{1/4}, \mu^{1/4})^{} \| e_{\lambda}\|_{L^2} \|e_{\mu}\|_{L^2}$$
on compact two-dimensional manifolds without boundary (this has been extended to higher dimensions \cite{BSSmulti, burq2}). A recent result of the third author \cite{stein}  (see also \cite{alex})
shows that one would generically (i.e., on typical manifolds in the presence of quantum chaos) expect $e_i(x)e_j(x)$ to be mainly supported at eigenfunctions having their eigenvalue close to $\max\left\{\lambda_i, \lambda_j\right\}$ and that deviation
from this phenomenon, as in the case of Fourier series on $\mathbb{T}$ for example, requires eigenfunctions to be strongly correlated at the wavelength in a precise sense.\\

In this paper, we ask the question on the numerical
rank of the space spanned by the pointwise products of eigenfunctions
$$ A_n = \mbox{span} \left\{ e_i(x) e_j(x): 1 \leq i,j \leq n\right\}.$$
This is a natural quantity for measuring the complexity of the products but also motivated by the density fitting approximation to the electron repulsion integral in the quantum chemistry literature. 
Given a set of eigenfunctions, the four-center two-electron repulsion integral
$$ (ij|kl) = \int_{\Omega \times \Omega}{\frac{e_i(x) e_j(x) e_k(y) e_l(y)}{|x-y|} dx dy}$$
is a central quantity in electronic structure theories. If we are
working with the first $n$ eigenfunctions
$(e_i)_{1 \leq i \leq n}$, then one has to evaluate $\mathcal{O}(n^4)$ integrals.

It has been empirically observed in the literature (see e.g.,
\cite{jian1} by the first author and Lexing Ying) that the space $A_n$
can in practice be very well approximated by another vector space
$B_n$ with $\mbox{dim}(B_n) \sim c \cdot n$, often referred as density
fitting in quantum chemistry literature.  This then drastically
reduces the number of integrals in need of evaluation to
$\mathcal{O}(n^2)$ and can be used for fast algorithms for electronic
structure calculations as in e.g., \cite{jian2}. Our result is
inspired by the empirical success of density fitting and gives a
mathematical justification.

\section{Main Results}
Motivated by the above, for $1\le i\le j\le n$, we are interested in
estimating products of eigenfunctions $e_ie_j$ by finite linear
combinations of eigenfunctions.  With this in mind, let for
$\nu\in {\mathbb N}$
\begin{equation}\label{1}
E_\nu f=\sum_{k=0}^\nu \langle f,e_k\rangle \, e_k
\end{equation}
denote the projection of $f\in L^2(M)$ onto the space, $B_\nu$, spanned by $\{e_k\}_{k=0}^\nu$, and
\begin{equation}\label{2}
\Rnu f=f-E_\nu f
\end{equation}
denote the ``remainder term'' for this projection.
Thus,
\begin{equation}\label{3}
\|f-E_\nu f\|^2_{L^2(M)}=\|\Rnu f\|_2^2=\sum_{k>\nu} |\langle f,e_k\rangle|^2.
\end{equation}

Our first result says that, if, as above, $i,j\le n$ then the ``eigenproduct'' $e_ie_j$ can be well approximated
by elements of $B_\nu$ if $\nu$ is not much larger than $n$.

\begin{theorem}\label{theorem1}   Fix $(M,g)$ as above.   Then there is a $\sigma=\sigma_d$ so that 
given $\kappa\in {\mathbb N}$ there is a uniform constant $C_\kappa$  such that if $i,j\le n<\nu$
and $n\in {\mathbb N}$ we have
\begin{equation}\label{4}
\|\Rnu (e_ie_j)\|_{L^2(M)}\le C_\kappa n^\sigma \, \bigl(n/\nu\bigr)^\kappa.
\end{equation}
Furthermore, there is a fixed $\sigma_\infty$ depending only on $d$ so that
\begin{equation}\label{5}
\|\Rnu (e_ie_j)\|_{L^\infty(M)}\le C_\kappa n^{\sigma_\infty} \, \bigl(n/\nu\bigr)^\kappa.
\end{equation}
Thus, given $\delta>0$ there is a constant $C_{M,\delta}$ so that if $\e\in (0,1)$ we have for $n\gg 1$
\begin{equation}\label{5'}
\|\Rnu (e_ie_j)\|_{L^\infty(M)} <\e, \quad \text{if } \, \,  \nu =\lceil C_{M,\delta} \, n^{1+\delta}\e^{-\delta} \rceil.
\end{equation}
\end{theorem}


Note that the results in Theorem~\ref{theorem1} trivially hold on the torus ${\mathbb T}^d \simeq (-\pi,\pi]^d$, since, in this case, $e_ie_j$ must be a trigonometric polynomial of
degree $\la_i+\la_j\le 2\la_n$ and $\Rnu$ annihilates trigonometric polynomials of degree $\la_\nu$.  Thus, $\Rnu(e_ie_j)=0$ in this case if $\nu$ is lager than a fixed
multiple of $n$.  The result is also trivial on the round sphere for similar reasons since the product of a spherical harmonic of degree $i$ with one of degree $j$ is a linear
combination of ones of degree $\le i+j$ whose expansion involves the Clebsch-Gordan coefficients.\\

On a general manifold there are no simple representations for eigenproducts and, in particular, it seems rare that the product of two eigenfunctions of frequency
$\le \la_n$ or less can be expressed by linear combinations of eigenfunctions of frequency $\la_\nu$ or less with $\nu$ being a fixed multiple of $n$.  On the other hand,
our result says that, in an approximate sense, one has a weaker version if one is willing to replace $\nu\approx n$ by $\nu\approx n^{1+\delta}$ for any $\delta>0$.
As we shall see, the proof of Theorem~\ref{theorem1} will also allow us to show that the same results hold for eigenproducts of any fixed length.

\begin{theorem}\label{theorem2}  Fix $\ell =2,3,\dots$.  Then there is a $\sigma_{d,\ell}<\infty$ so that if $j_1, \dots, j_\ell \le n < \nu$, with $n$, $\nu\in {\mathbb N}$, we have
\begin{equation}\label{6}
\|\Rnu(e_{j_1}\cdots e_{j_\ell})\|_{L^\infty(M)} \le C_{M,\kappa,\ell} \, \, n^{\sigma_{d,\ell}} \, \bigl(n/\nu\bigr)^\kappa, \quad \kappa=1,2,\dots,
\end{equation}
with the uniform constant  $C_{M,\kappa,\ell} $ depending only on $(M,g)$, $\kappa$ and $\ell$.
Thus, given $\delta>0$ we have for $\e \in (0,1)$
\begin{equation}\label{6'}
\|\Rnu(e_{j_1}\cdots e_{j_\ell})\|_{L^\infty(M)} <\e, \quad \text{if } \nu =\lceil C_{M, \ell,\delta} \, n^{1+\delta}\e^{-\delta} \rceil.
\end{equation}
\end{theorem}

To prove these results we shall use a classical Sobolev embedding; for products with more terms, we use a variation that comes from a combination of Leibniz's rule together with estimates established by the second author  \cite{sogge88}.

\begin{lemma}\label{lemma}
For $\sigma\in \R$, let  $\|f\|_{H^\sigma(M)}=\|(I-\Delta_g)^{\sigma/2}f\|_{L^2(M)}$ denote the norm for the Sobolev space of order $\sigma$ on $M$.
Then
\begin{equation}\label{8}
\|f\|_{L^\infty(M)}\lesssim \|f\|_{H^\sigma(M)}, \quad \text{if } \, \, \sigma>d/2.
\end{equation}
Also, if $\ell, \, n,\, \mu\in {\mathbb N}$ and if $1\le j_1,\dots, j_\ell \le n$,
\begin{equation}\label{9}
\bigl\| \, \prod_{k=1}^\ell e_{j_k}\, \bigr\|_{H^\mu(M)} \le C_{\ell,\mu,M} \, \la_n^{\mu+\ell \sigma(2\ell,d)},
\end{equation}
if, for $p>2$ we set
\begin{equation}\label{10}
\sigma(p,d)=\max\bigl\{ \tfrac{d-1}2\bigl(\tfrac12-\tfrac1p\bigr), \, d\bigl(\tfrac12-\tfrac1p)-\tfrac12 \, \bigr\}.
\end{equation}
\end{lemma}

These bounds arise naturally from estimates established by the second author  \cite{sogge88}
saying that if $p>2$ then for $j\ge 1$ we have
\begin{equation}\label{9'}
\|e_j\|_{L^p(M)}\lesssim \la_j^{\sigma(p,d)},
\end{equation}
with $\sigma(p,d)$ being as in \eqref{10}.

\section{Proofs}
\subsection{Proof of Theorem 1 and Theorem 2}
\begin{proof}
To prove the $L^2$-estimate \eqref{4} we note that if $\mu\in {\mathbb N}$, we have
\begin{align*}
\|\Rnu h\|^2_{L^2}=\sum_{k>\nu} \bigl|\langle h,e_k\rangle\bigr|^2 &\le \la_\nu^{-2\mu}\sum_{k>\nu}\la_k^{2\mu}\, \bigl|\langle h,e_k\rangle\bigr|^2
\\
&\le \la_\nu^{-2\mu}\|h\|^2_{H^\mu}.
\end{align*}
If we take $h=e_ie_j$ and use this along with \eqref{9} we conclude that
$$\|\Rnu(e_ie_j)\|_{L^2}\lesssim \la_\nu^{-\mu} \,  \la_n^{2\sigma(4,d)+\mu}.$$
Since, by the Weyl formula, $n\approx \la_n^d$ and $\nu\approx \la_\nu^d$,  this inequality implies that
$$\|\Rnu (e_ie_j)\|_{L^2(M)}\le C \, \bigl(n/\nu\bigr)^{\mu/d} \, n^{\frac2d \sigma(4,d)}.
$$
As $\mu\in {\mathbb N}$ is arbitrary, this 
 of course yields \eqref{4} with $\sigma$ there being $\tfrac2d \sigma(4,d)$.
To prove the sup-norm bounds assume that $\mu\in {\mathbb N}$ is larger than $(d+1)/2$.  Then, by \eqref{8} and the above argument we have
\begin{align*}
\|\Rnu (e_ie_j)\|^2_\infty &\lesssim \|\Rnu(e_ie_j)\|_{H^{(d+1)/2}}^2 \\
&\le \la_\nu^{d+1-2\mu} \|e_ie_j\|_{H^\mu}^2\\
&\lesssim  \la_\nu^{d+1-2\mu} \la_n^{4\sigma(4,d)+2\mu}.
\end{align*}
From this, we obtain
\begin{align*}
\|\Rnu (e_ie_j)\|_\infty &\lesssim \la_n^{2\sigma(4,d)+(d+1)/2} \bigl(\la_n/\la_\nu\bigr)^{\mu-(d+1)/2} \\
&\approx n^{\frac2d\sigma(4,d)+\frac{d+1}{2d}} \, \bigl(n/\nu\bigr)^{[\mu-(d+1)/2]/d},
\end{align*}
which yields \eqref{5} with 
$$\sigma_\infty =\tfrac2d \sigma(4,d)+\tfrac{d+1}{2d}.$$
Also, \eqref{6'} is a trivial consequence of \eqref{6}.
This argument clearly also gives the approximation bounds \eqref{6} for eigenproducts of length $\ell$.  Indeed by \eqref{8}, if 
$\mu\in {\mathbb N}$ is larger than $(d+1)/2$,
$$\|\Rnu(e_{j_1}\cdots e_{j_\ell})\|_\infty \le \la^{(d+1)/2-\mu}_\nu \bigl\| \, \prod_{k=1}^\ell e_{j_k}\, \bigr\|_{H^\mu}
\lesssim  \la^{(d+1)/2-\mu}_\nu \la_n^{\mu+\ell \sigma(2\ell,d)},
$$
which yields \eqref{6} with $$\sigma_{d,\ell}=\tfrac\ell{d}\sigma(2\ell,d)+\tfrac{d+1}{2d}.$$
\end{proof}

\subsection{Proof of Lemma~\ref{lemma}}
\begin{proof}
To prove \eqref{9} we first recall some basic facts about Sobolev spaces on manifolds.  See \cite{SFIOII} for more detals.
First, if $1=\sum_{j=1}^N \varphi_j$ is a fixed smooth partition of unity with
$$\text{supp }\varphi_j  \Subset \Omega_j,$$
where $\Omega_j\subset M$ is a coordinate patch, we have for fixed $\mu\in {\mathbb N}$
$$\|f\|_{H^\mu(M)} \approx \sum_{j=1}^N \sum_{|\alpha|\le \mu} \bigl\| \partial^\alpha (\varphi_j f)\bigr\|_{L^2(\Rn)}.$$
Here, 
the $L^2$-norms are taken with respect to our local coordinates.
 $\| \prod_{k=1}^\ell e_{j_k}\|_{H^\mu(M)}$ is dominated by a finite sum of terms of the form
$$\|\partial^\alpha (\varphi \cdot \prod_{k=1}^\ell e_{j_k})\|_{L^2},$$
where $\varphi=\varphi_j$ for some $j=1,\dots,N$ and $|\alpha|\le \mu$.
By Leibniz's rule, we can thus dominate the left side of \eqref{9} by a finite sum of terms of the form
$$\bigl\| \prod_{k=1}^\ell L_k e_{j_k}\bigr\|_{L^2(M)},$$
where $L_k: \, C^\infty(M)\to C^\infty(M)$ are differential operators with smooth coefficients of order $m_k$ with
\begin{equation}\label{12}
m_1+\cdots+m_\ell \le \mu.
\end{equation}
As a result, by H\"older's inequality, $\|\prod_{k=1}^\ell e_{j_k}\|_{H^\mu(M)}$ is majorized by a finite sum of terms of the form
\begin{equation}\label{13}
\prod_{k=1}^\ell \|L_ke_{j_k}\|_{L^{2\ell}(M)},
\end{equation}
where the $L_k$ are as above.
Since $L_k$ is a differential operator of order $m_k$, for any $1<p<\infty$, standard $L^p$ elliptic regularity estimates give
$$\|L_k h\|_{L^p(M)} \lesssim \|(I-\Delta_g)^{m_k/2}h\|_{L^p(M)}.$$
Thus, since $(I-\Delta_g)^{m_k/2}e_{j_k}=(1+\la_{j_k}^2)^{m_k/2}e_{j_k}$ and $0<\la_1\le \la_{j_k}\le \la_n$
$$\|L_k e_{j_k}\|_{L^{2\ell}(M)} \lesssim \la_n^{m_k}\|e_{j_k}\|_{L^{2\ell}(M)} \lesssim \la_n^{m_k+\sigma(2\ell,d)},$$
using \eqref{9'} in the last inequality.
Since \eqref{13} has $\ell$ factors and \eqref{12} is valid we obtain \eqref{9} from this, which finishes the proof of Lemma~\ref{lemma}.
\end{proof}

\subsection{Some Remarks}  We just showed that, if $\sigma(p,d)$ is as in \eqref{10} then if $n, \nu \in {\mathbb N}$ and
$j_1, \dots, j_\ell \le n$ and $n<\nu$, then we have
\begin{equation}\label{14}
\bigl\| \Rnu( \prod_{k=1}^\ell e_{j_k})\bigr\|_{L^2(M)}\le C_\kappa \la_n^{\ell \sigma(2l,d)} \, \bigl(\la_n/\la_\nu \bigr)^\kappa
\approx n^{\ell \sigma(2\ell,d)/d} \, (n/\nu)^{\kappa/d},
\end{equation}
for each $\kappa=1,2,3,\dots$.  Since $\|\Rnu\|_{L^2\to L^2}=1$, by H\"older's inquality and \eqref{9'} we clearly have
\begin{equation}\label{15}
\bigl\| \Rnu( \prod_{k=1}^\ell e_{j_k})\bigr\|_{L^2(M)}\le C \la_n^{\ell \sigma(2l,d)}, \quad \text{if } \, \, \nu\le n.
\end{equation}
Thus, if $\kappa\in {\mathbb N}$ is fixed, we have the uniform bounds
\begin{align}\label{16}
\bigl\| \Rnu( \prod_{k=1}^\ell e_{j_k})\bigr\|_{L^2(M)}
&\le C_\kappa \la_n^{\ell \sigma(2l,d)}\bigl(1+\la_\nu/\la_n\bigr)^{-\kappa} \\
&\approx n^{\ell \sigma(2\ell,d)/d} \bigl(1+\nu/n\bigr)^{-\kappa/d}.
\end{align}

Inequality \eqref{16} is saturated on the sphere $S^d$
for $\nu$ smaller than a fixed multiple of $n$  if $j_1=\cdots = j_\ell =n$.  Specifically, if  $\ell \ge 3$ or if $\ell =2$ and $d\ge 3$ one can check that zonal
functions saturate the bound.  For the remaining case where $\ell = d=2$, the highest weight spherical harmonics saturate the bounds if $\la_\nu+1\le 2\la_n$.  On the other
hand, as we pointed out before, the left side of \eqref{16} is zero in this case if $\la_\nu+1>2\la_n$.


The proof of Theorems~\ref{theorem1}--\ref{theorem2} shows that if $1\le j_k\le n< \nu$ and $\kappa \in {\mathbb N}$, then
$$\bigl\| \Rnu(\prod_{k=1}^\ell e_{j_k})\bigr\|_{L^2(M)} \le C_\kappa \bigl(\la_\nu/\la_n\bigr)^{-\kappa}\prod_{k=1}^\ell \|e_{j_k}\|_{L^{2\ell}(M)}.$$
We then used the universal $L^p$-bounds \eqref{9'} of one of us to control the right side and prove our results.
Substantially improved eigenfunction estimates would of course lead to better bounds for the $L^2$-approximation of products of eigenfunctions.  For instance, the ``random wave model''
 of Berry~\cite{Berry}  predicts that eigenfunctions on Riemannian manifolds with chaotic geodesic flow should have $\|e_\la\|_{L^p}=O(1)$ for $p<\infty$.  (See e.g., \cite{ZelSurv} for
 more details.)  If this optimistic conjecture were valid we would have, for any fixed $\ell =2,3,\dots$, 
 $$\bigl\| \Rnu (\prod_{k=1}^\ell e_{j_k})\bigr\|_2 <\e \quad \text{if }\, \, \nu>C_\delta \, n \, \e^{-\delta} \, \, \, \text{and } \, \, 1\le j_k\le n,$$
 with $\delta>0$ being arbitrary.  We should note, though, that there has been much work 
 on obtaining improved $L^p$-estimates for eigenfunctions over the
 last forty years and only logarithmic improvements over the universal
 bounds \eqref{9'} have been obtained assuming, say, that one has
 negative curvature.  See, e.g., \cite{Berard, SBLog,
   BlairSoggeKaknik, HassellTacyNonpos, SoCrit, SoggeZelditchMaximal,
   SoggeZelditchL4}.

We also would like to point out that the argument that we have employed to obtain the approximation bounds in Theorem~\ref{theorem1} yield higher order Sobolev estimates as well.
Indeed, if $\ell=2,3,\dots$ and $\mu \in {\mathbb N}$ are fixed, the proof of \eqref{5'}  shows that if, $1\le j_k\le n$, $k=1,\dots, \ell$, 
\begin{equation}\label{21}
\bigl\| \Rnu(\prod_{k=1}^\ell e_{j_k})\bigr\|_{H^\mu}<\e \quad
\text{for } \, \e \in (0,1) \quad \text{if } \, \, \nu=O( n^{1+\delta} \, \e^{-\delta}),
\end{equation}
where $\delta>0$ can be chosen to be as small as we like.  This in turn allows us to control the $C^m$ norms of $\Rnu(e_{j_1}\cdots e_{j_\ell)}  $ for any $m\in {\mathbb N}$ for such $\nu$.


\section{Approximation in $H^{-1}$}

Motivated by the discussion of the four-center two-electron repulsion integral, which is a somewhat better behaved quantity due to the smoothing effects of the potential, it is
also natural to look for an approximation result in a function space that captures the smoothing effect of the potential.  Since a multiple of  $|x-y|^{-1}$ is the fundamental solution of the
Laplacian in ${\mathbb R}^3$, the appropriate physically relevant problems involve the Sobolev space $H^{-1}$ equipped with the norm defined by
\begin{equation}\label{3.0}
\|f\|_{H^{-1}}^2=\bigl\| (I-\Delta_g)^{-1/2}f\|_{L^2}^2=\sum (1+\la_k^2)^{-1} |\langle f,e_k\rangle|^2.
\end{equation}

Our $H^{-1}$ approximation result then is the following.

\begin{theorem}\label{theoremh}  
Let $(M,g)$ be a compact Riemannian manifold of dimension $2\le d\le 4$.  Then if $1\le i,j\le n<\nu$ and $\e \in (0,1)$ we have
\begin{equation}\label{3.1}
\bigl\| \Rnu(e_ie_j)\bigr\|_{H^{-1}}<\e, \qquad \text{if } \, \,  \nu=O(n^{\mu(d)}\e^{-d}),
\end{equation}
where
\begin{equation}\label{3.2}
\mu(d)=
\begin{cases}
1/4, \quad \text{if } \, d=2
\\
1/2,  \quad \text{if } \, d=3
\\
1, \,   \,  \,  \, \,  \quad \text{if } \,  d=4.
\end{cases}
\end{equation}
\end{theorem}
Thus, our results only involve sublinear growth of $\nu$ in terms of $n$ for $d\le 3$ and linear growth for $d=4$, which matches up nicely with the trivial cases of ${\mathbb T}^d$ and
$S^d$ for these dimensions.
We only stated things in dimensions $d\le 4$ since in higher dimensions we cannot get improved $H^{-1}$ bounds compared to the $L^2$ bounds in Theorem~\ref{theorem1}.

\begin{proof} Under the above hypothesis, we claim that
\begin{equation}\label{3.3}
\bigl\|\Rnu(e_ie_j)\bigr\|_{H^{-1}}\le C\la^{-1}_\nu \la_n^{2\sigma(4,d)} \approx \nu^{-1/d}n^{2\sigma(4,d)/d},
\end{equation}
where, as in \eqref{10},
\begin{equation}\label{3.4}
\sigma(4,d)=\begin{cases}
\tfrac18, \qquad \, \text{if } \, d=2
\\
\tfrac{d-2}4, \,  \quad \text{if } \, d\ge 3.
\end{cases}
\end{equation}
Since 
$$\nu^{-1/d}n^{2\sigma(4,d)/d}<\e \,  \iff  \, \nu>\e^{-d}n^{\mu(d)},$$
 where $\mu(d)$, as in \eqref{3.2}, equals $2\sigma(4,d)$, we conclude that
in order to obtain \eqref{3.1},
we just need to prove \eqref{3.3}.
To prove this we  use H\"older's inequality and \eqref{9'} to get
\begin{align*}
\bigl\|\Rnu(e_ie_j)\bigr\|_{H^{-1}}^2
&\le \la^{-2}_\nu\sum_{k>\nu}  |\langle e_ie_j,e_k\rangle|^2 \\
&\le \la^{-2}_\nu \|e_i e_j\|_2^2 \le \la^{-2}_\nu \|e_i\|_4^2 \, \|e_j\|_4^2 \lesssim \la^{-2}_\nu \la_n^{4\sigma(4,d)},
\end{align*}
as desired.
\end{proof}

\section{Results for manifolds with boundary and domains}  

Let us conclude by extending our results to
eigenfunctions of the Dirichlet Laplacian on $M$, where $M$ is either a $d$-dimensional relatively compact
Riemannian manifold or a bounded domain in ${\mathbb R}^d$ and where the boundary, $\partial M$,
of $M$ is smooth.  Thus, we shall consider an $L^2$-normalized basis of eigenfunctions $\{e_j\}_{j=1}^\infty$ satisfying
$$-\Delta_g e_j(x)=\la_j^2e_j(x), \, \, x\in M, \quad \, \text{and } \, e_j|_{\partial M}=0.$$
As before, we shall assume that the frequencies of the eigenfunctions are arranged in increasing order, i.e.,  $0<\la_1\le \la_2\le \la_3\le\dots$.

If we then define $\Rnu$ as in \eqref{2} we have the following.

\begin{theorem}\label{theoremd}
Let $\ell =2,3,\dots$.  Assume that $j_k\le n$, $k=1,\dots, \ell$.  Then given $\delta>0$ we have
for $\e\in (0,1)$
\begin{equation}\label{4.1}
\bigl\| \Rnu(\prod_{k=1}^\ell e_{j_k})\bigr\|_{L^\infty} <\e, \qquad
\text{if } \, \, \, \nu=O(n^{1+\delta}\e^{-\delta}).
\end{equation}
Also, if $2\le d\le 4$ and
$\|f\|_{H^{-1}}$ is as in \eqref{3.0} and $i,j\le n\in {\mathbb N}$
\begin{equation}\label{4.2}
\bigl\|\Rnu(e_ie_j)\bigr\|_{H^{-1}}<\e, \qquad \text{if } \, \, \, \nu=O(n^{\mu(d)}\e^{-d}),
\end{equation}
where $\mu(2)=1/3$, $\mu(3)=2/3$ and $\mu(4)=1$.
\end{theorem}

To prove these results one uses the eigenfunction estimates in \cite{SmSoActa} which say that bounds
of the form \eqref{9'} are valid, but where the exponent $\sigma(p,d)$ may be larger than the one
in \eqref{10} for a certain range of $p$ depending on $d$.  Our earlier argument only used the fact that 
this exponent was finite in order to get bounds of the form \eqref{4.1} and the same applies here.  In fact
by using such a bound along with elliptic regularity estimates (see \cite[\S 9.6]{Gilbarg}), one obtains
\eqref{4.1}. 

By the same argument we can also obtain bounds of the form \eqref{21} in this setting.

To prove \eqref{4.2} one, as before, only needs to use $L^4$ eigenfunction estimates.  Specifically, by the
arguments from the preceding section, we obtain \eqref{4.2} just by using the fact that, for Dirichlet eigenfunctions, we have \eqref{9'} for
$p=4$ with
\begin{equation*}
\sigma(4,d)=
\begin{cases}
1/6, \quad \text{if } d=2
\\
1/3, \quad \text{if } d=3
\\
1/2, \quad \text{if } d=4.
\end{cases}
\end{equation*}
These bounds were obtained in \cite{SmSoActa}.

We remark that if $M$ has geodesically concave boundary, the results in \cite{SmSoCrit} say that we have the more favorable bounds where,
as in the boundaryless case, $\sigma(4,d)$ is given by \eqref{3.4}  and so, in this case, we can
recover the bounds in \eqref{3.1}--\eqref{3.2}.

\bibliographystyle{acm}

\bibliography{refs}

\providecommand{\MR}[1]{}
\begin{thebibliography}{10}

\bibitem{Berard}
{\sc B\'erard, P.~H.}
\newblock On the wave equation on a compact {R}iemannian manifold without
  conjugate points.
\newblock {\em Math. Z. 155}, 3 (1977), 249--276.

\bibitem{bern}
{\sc Bernstein, J., and Reznikoff, A.}
\newblock Analytic continuation of representations and estimates of automorphic
  forms.
\newblock {\em Ann. Math. 150\/} (1999), 329--352.

\bibitem{Berry}
{\sc Berry, M.~V.}
\newblock Regular and irregular semiclassical wavefunctions.
\newblock {\em J. Phys. A 10}, 12 (1977), 2083--2091.

\bibitem{BSSmulti}
{\sc Blair, M.~D., Smith, H.~F., and Sogge, C.~D.}
\newblock {On multilinear spectral cluster estimates for manifolds with
  boundary}.
\newblock {\em Math. Res. Lett. 15}, 3 (2008), 419--426.

\bibitem{SBLog}
{\sc Blair, M.~D., and Sogge, C.~D.}
\newblock Logarithmic improvements in {$L^p$} bounds for eigenfunctions at the
  critical exponent in the presence of nonpositive curvature.
\newblock arXiv:1706.06704.

\bibitem{BlairSoggeKaknik}
{\sc Blair, M.~D., and Sogge, C.~D.}
\newblock {On {K}akeya-{N}ikodym averages, {$L^p$}-norms and lower bounds for
  nodal sets of eigenfunctions in higher dimensions}.
\newblock {\em J. Eur. Math. Soc. (JEMS) 17}, 10 (2015), 2513--2543.

\bibitem{burq2}
{\sc Burq, N., G\'{e}rard, P., and Tzvetkov, N.}
\newblock Multilinear estimates for the laplace spectral projectors on compact
  manifolds.
\newblock {\em C. R. Math. Acad. Sci. Paris 338\/} (2004), 359--364.

\bibitem{burq}
{\sc Burq, N., G\'{e}rard, P., and Tzvetkov, N.}
\newblock Bilinear eigenfunction estimates and the nonlinear {S}chr\"odinger
  equation on surfaces.
\newblock {\em Invent. Math. 159\/} (2005).

\bibitem{alex}
{\sc Cloninger, A., and Steinerberger, S.}
\newblock On the dual geometry of {L}aplacian eigenfunctions, 2018.
\newblock preprint, arXiv:1804.09816.

\bibitem{Gilbarg}
{\sc Gilbarg, D., and Trudinger, N.~S.}
\newblock {\em Elliptic partial differential equations of second order},
  second~ed., vol.~224 of {\em Grundlehren der Mathematischen Wissenschaften
  [Fundamental Principles of Mathematical Sciences]}.
\newblock Springer-Verlag, Berlin, 1983.

\bibitem{HassellTacyNonpos}
{\sc Hassell, A., and Tacy, M.}
\newblock {Improvement of eigenfunction estimates on manifolds of nonpositive
  curvature}.
\newblock {\em Forum Mathematicum 27}, 3 (2015), 1435--1451.

\bibitem{krotz}
{\sc Kr\"otz, B., and Stanton, R.}
\newblock Holomorphic extension of representations: ({I}) automorphic
  functions.
\newblock {\em Ann. Math. 159\/} (2004), 641--724.

\bibitem{jian2}
{\sc Lu, J., and Thicke, K.}
\newblock Cubic scaling algorithms for rpa correlation using interpolative
  separable density fitting.
\newblock {\em J. Comput. Phys. 351\/} (2017), 187--202.

\bibitem{jian1}
{\sc Lu, J., and Ying, L.}
\newblock Compression of the electron repulsion integral tensor in tensor
  hypercontraction format with cubic scaling cost.
\newblock {\em J. Comput. Phys. 302\/} (2015), 329--335.

\bibitem{Sarnak}
{\sc Sarnak, P.}
\newblock Integrals of products of eigenfunctions.
\newblock {\em IMRN 6\/} (1994), 251--260.

\bibitem{SmSoCrit}
{\sc Smith, H.~F., and Sogge, C.~D.}
\newblock {On the critical semilinear wave equation outside convex obstacles}.
\newblock {\em J. Amer. Math. Soc. 8}, 4 (1995), 879--916.

\bibitem{SmSoActa}
{\sc Smith, H.~F., and Sogge, C.~D.}
\newblock {On the {$L^p$} norm of spectral clusters for compact manifolds with
  boundary}.
\newblock {\em Acta Math. 198}, 1 (2007), 107--153.

\bibitem{sogge88}
{\sc Sogge, C.~D.}
\newblock {Concerning the {$L^p$} norm of spectral clusters for second-order
  elliptic operators on compact manifolds}.
\newblock {\em J. Funct. Anal. 77}, 1 (1988), 123--138.

\bibitem{SFIOII}
{\sc Sogge, C.~D.}
\newblock {\em Fourier integrals in classical analysis}, second~ed., vol.~210
  of {\em Cambridge Tracts in Mathematics}.
\newblock Cambridge University Press, Cambridge, 2017.

\bibitem{SoCrit}
{\sc Sogge, C.~D.}
\newblock Improved critical eigenfunction estimates on manifolds of nonpositive
  curvature.
\newblock {\em Math. Res. Lett. 24}, 2 (2017), 549--570.

\bibitem{SoggeZelditchMaximal}
{\sc Sogge, C.~D., and Zelditch, S.}
\newblock {Riemannian manifolds with maximal eigenfunction growth}.
\newblock {\em Duke Math. J. 114}, 3 (2002), 387--437.

\bibitem{SoggeZelditchL4}
{\sc Sogge, C.~D., and Zelditch, S.}
\newblock {On eigenfunction restriction estimates and {$L^4$}-bounds for
  compact surfaces with nonpositive curvature}.
\newblock In {\em {Advances in analysis: the legacy of {E}lias {M}. {S}tein}},
  vol.~50 of {\em {Princeton Math. Ser.}} Princeton Univ. Press, Princeton, NJ,
  2014, pp.~447--461.

\bibitem{stein}
{\sc Steinerberger, S.}
\newblock On the spectral resolution of products of {L}aplacian eigenfunctions.
\newblock {\em J. Spectral Theory\/} (in press).

\bibitem{ZelSurv}
{\sc Zelditch, S.}
\newblock {Quantum ergodicity and mixing of eigenfunctions}.
\newblock In {\em {Elsevier Encyclopedia of Math. Phys.}}, vol.~1. Elsevier,
  2006, pp.~183--196.

\end{thebibliography}

%

\end{document}